\documentstyle[12pt]{article}
\input amssym.def
\input amssym.tex 

\newenvironment{demo}[1]%
{\vskip-\lastskip\medskip
  \noindent
  {\em #1.}\enspace
  }%
{\qed\par\medskip
  }

\newcommand{\qed}{
  \strut\hfill
  \mbox{$\Box$}
  }

{\vskip-\lastskip
  \begin{trivlist}
    \item[]
      {\bf Axiom #1 }
      }%
{\end{trivlist}
  }

\newtheorem{theorem}{Theorem}[section]

\newtheorem{corollary}{Corollary}[section]

\newtheorem{lemma}{Lemma}[section]

\newtheorem{remark}{Remark}[section]

\newtheorem{proposition}{Proposition}[section]

\newcommand{\bg}{
  \beta\gamma
  }

\newcommand{\bngw}{
  :\gamma (w)\partial^n \beta (w):
  }

\newcommand{\bNg}{
  :\gamma (z)\partial^N \beta (z):
  }

\newcommand{\bIg}{
  :\gamma (z) \partial^i \beta (z):
  }

\newcommand{\bi}{
  {
    n \choose k
    } 
  }

\newcommand{\D}{
  \cal D
  }

\newcommand{\F}{
  {\cal F}
  }

\newcommand{\Fl}{
  {\cal F}^l
  }

\newcommand{\Fbar}{
  {\overline{\cal F} }^l
  }

\newcommand{\Hil}{
  {\cal H}_{i(s+l)}
  }
  
\newcommand{\hD}{
  { \widehat{\cal D} }
  }

\newcommand{\hf}{
  \frac12
  }

\newcommand{\Ho}{
  {\cal H}_0
  }

\newcommand{\Ltw}{
   {\cal L}_c (t, w)
  }

\newcommand{\M}{
   {\cal M}_s
  }

\newcommand{\Moo}{
   {\cal M} (0, 0)
  }

\newcommand{\Mtw}{
   {\cal M}_c (t, w)
  }

\newcommand{\Ml}{
   {\cal M}_s^l
  }

\newcommand{\N}{
  (n)z^{-n + s}
  }

\newcommand{\NN}{
  (n)z^{-n -s -1}
  }

\newcommand{\ps}{
  \psi (z)
  }

\newcommand{\ph}{
  \phi (z)
  }

\newcommand{\SUM}{
  \sum_{n\in \Bbb Z}
  }

\newcommand{\Tr}{
  \mbox{Tr}
  }

\newcommand{\UW}{
  {\cal U} ({\cal W}_{3, c})
  }

\newcommand{\vac}{
  | 0 \rangle
  }

\newcommand{\vacoo}{
  | 0,0 \rangle
  }

\newcommand{\vacs}{
  | s \rangle
  }

\newcommand{\vacbc}{
  | bc \rangle
  }

\newcommand{\vactw}{
  | t, w \rangle
  }

\newcommand{\vacuum}{
  | \alpha \rangle
  }

\newcommand{\VM}{
  {\cal VW}_{3, c}
  }

\newcommand{\W}{
  { \cal W }_{1+\infty}
  }

\newcommand{\Wc}{
  {\cal W}_{{1 + \infty}, c}
}

\newcommand{\Wth}{
  {\cal W}_3
  }

\newcommand{\Wtho}{
  {\cal W}_{3,0}
  }

\newcommand{\Wthp}{
  {\cal W}_{3, +}
  }

\newcommand{\Wthpm}{
  {\cal W}_{3, \pm}
  }

\newcommand{\Wthtwo}{
  {\cal W}_{3, -2}
  }

\newcommand{\WGth}{
  {\cal{W} } (gl_3)
  }

\newcommand{\WGN}{
  {\cal{W} } (gl_N)
  }

\newcommand{\WinfN}{
  {\cal W}_{{1 + \infty}, N}
  }

\newcommand{\Winfone}{
  {\cal W}_{{1 + \infty}, -1}
  }

\newcommand{\WminusN}{
  {\cal W}_{{1 + \infty}, -N}
  }

\newcommand{\WN}{
  {\cal W}_N
  }

\newcommand{\WSN}{
  {\cal{W} } (sl_N)
  }

\newcommand{\Z}{
  \Bbb Z
  }

\begin{document}
\title{
  $\W\/$ algebra, $\Wth$ algebra, and Friedan-Martinec-Shenker\\
bosonization
  }
\author{
  %
  Weiqiang Wang\\
\\{\small Max-Planck Institut fur Mathematik, 
53225 Bonn, Germany}\\
{\small E-mail: wqwang@mpim-bonn.mpg.de}
}
\date{}
\maketitle
\begin{abstract}
  We show that the vertex algebra $\W$ with central
charge $-1$ is isomorphic to a tensor product
of the simple $\Wth$ algebra with central charge 
$-2$ and a Heisenberg vertex algebra generated by
a free bosonic field. We construct a family of irreducible
modules of the $\Wth$ algebra with central charge 
$-2$ in terms of free fields and calculate the
full character formulas
of these modules with respect to the full
Cartan subalgebra of the $\Wth$ algebra. 
\end{abstract}

\setcounter{section}{-1}
\section{Introduction}

In search of classification of conformal field theories, one 
is lead to study $\cal W\/$ algebras which
are extended chiral algebras (vertex algebras
or vertex operator algebras
in mathematical terminology) containing
Virasoro algebra as a subalgebra. 
Since the first attempt
was made by Zamolodchikov \cite{Z} there has been
much further study of $\cal W$ algebras
(see the review paper \cite{BS}
and references therein)\footnote{We note a less-known
fact that $\WN$ algebras were constructed in \cite{F2}
for the particular central charge $c = N-1$}. A particularly interesting
example of $\cal W$ algebra, the so-called
$\W\/$ algebra \cite{PRS}, appears to be a universal one
among various $\cal W\/$ infinite algebras
in the $ N \rightarrow \infty $
limit of $\WSN$ algebras, see e.g. \cite{Ba, BK, PRS, O}. 
The $\WSN$ algebras are often referred to as $\WN$ algebras in literature.
In mathematics, $\W$ is known as 
the universal central extension $\hD$ of the
Lie algebra $\D$ of differential operators on the
circle. The first systematic study of 
the representation theory of the Lie algebra
$\hD\/$ was undertaken by Kac and Radul 
in their seminal paper \cite{KR1} and
there have been many
further development \cite {M, FKRW, AFMO, KR2, W1}
since then, just to name some.

In \cite{FKRW}, Lie algebra
$\hD$ and its representation theory are
studied in the framework of vertex algebras \cite{B, FLM, DL, K2, LZ2}.
It turns out that the irreducible vacuum $\hD\/$-module
with central charge $c$ admits a canonical
vertex algebra structure, with infinitely
many generating fields of conformal weights $2,3,4,\cdots$, 
which we will denote by $\Wc$.
The case when the central charge is non-integral is 
not difficult to understand. The case
when the central charge is a positive integer 
was studied in detail in \cite {FKRW}. 
The vertex algebra $ \WinfN \/ $ with a positive 
integral central charge $N\/$ has redundant symmetries,
namely only the first $N$ generating fields are independent.
More precisely $ \WinfN \/ $ is shown to be isomorphic to 
a $\cal W\/$-algebra $\WGN$ with central
charge $N\/$ and the irreducible modules 
of $\WinfN$ are classified \cite{FKRW}. 

In this paper we will take the first step to clarify the 
connection between vertex algebra $\WminusN$ and some
other $\cal W$--algebra with finitely many generating fields. 
We prove that the vertex algebra $\Winfone$ is isomorphic to
a $\WGth$ algebra, which is a tensor product
of the simple $\Wth$ algebra with central charge $-2$ 
(denoted below by $\Wthtwo$) and 
a Heisenberg vertex algebra generated by a free bosonic field.
We will construct explicitly a number of modules of the
$\Wthtwo$ algebra parametrized by integers
in terms of free fields. 
We prove the irreducibility of these modules. 
As a by-product, we obtain full character formulas
for these representations. To our best knowledge, these
seem to be the first known full character formula of 
any non-trivial module of the $\Wth$ algebra 
with any non-generic central charge. We mention a curious
fact that a generating function of counting covers of
an elliptic curve \cite{Di} appears to be closely related to
our character formulas and admits very interesting
modular invariance properties \cite{KZ}.

The difficulties appearing in the negative integral central
charge case in contrast to the positive integral central
charge case are roughly the following: In both cases 
we have free field realizations. In the case
of positive integral central charge case we need
$bc$ fields which are free fermions 
while in the negative integral
central charge case we need $\beta\gamma$ fields which
are free bosonic ghosts. The structure of $\WN$ algebra
in the realization of $\WinfN$ in terms of
$bc$ fields can be identified 
relatively easily due to the very fact
that the structure of the basic
representation of the affine Kac-Moody algebra 
$\widehat{sl}_N\/$ is well understood \cite{K1}. However
structures of representations of 
affine algebras with negative integral central charges
are far from being clear. 

One of the main technique we use in 
relating the $\W$ algebra with central charge $-1$
to $\Wth$ algebra with central charge $-2$
is the bosonization of $\bg$ fields \cite{FMS}.
A similar construction was also given by
Kac and van de Leur and used by them for a construction
of a super KP hierarchy \cite{KV1, KV2}. More detailed
structures in the bosonization of
$\bg$ fields are further worked out in \cite{FF} and used 
for the computation of semi-infinite cohomology of
the Virasoro algebra with coefficient in the module of its adjoint
semi-infinite symmetric powers.
It is well known that $\bg\/$ fields 
are fundamental ingredients in superstring theory \cite{FMS},
in realizations of level $-1$ representations
of classical affine algebras
\cite{FeF} and in the calculation of BRST
cohomology of super-Virasoro algebras \cite{LZ1}.
They are also closely related to the logarithmic
conformal field theories which recently attract much
attention from physicists, see e.g. \cite{F, Ka, GK}. 
We hope our results may shed some lights on these subjects.

Let us explain in more detail. It is well known \cite{M} that the
Fock space $\M$ of the $\bg$ fields as a module over $\hD$
can be decomposed into a direct sum of the 
modules $\Ml$ parametrized by the $\bg$--charge number $l$.
Recall \cite{FMS} that the $\bg$ fields are expressed in terms of two
scalar fields $\ps$ and $\ph$. So the space $\M$ can be 
identified with some subspace of the Fock space of
the Heisenberg algebra of the two scalar fields $\ps$ and $\ph$.
Indeed one can identify $\Ml$ as
$\Fbar \bigotimes \Hil, \, i = \sqrt{-1}$ (cf. e.g. \cite{FF}), where $\Fbar $
is a certain subspace of the Fock space of the Heisenberg
algebra generated by the Fourier components of
field $\ps$ while $\Hil$ is the Fock
space of the Heisenberg algebra of the field $\ph$.

$\Winfone$ acts on $\Ml$ by means of fields
$$ J^i (z) =  \bIg+ \frac{1}{i+1} s(s-1)\cdots (s-i)z^{-i-1},
i \in \Z_{+}. $$
For the sake of simplicity, the reader may 
understand the main results of this paper by taking $ s =0$
throughout this paper.
By the celebrated boson-fermion correspondence, we have a
pair of fermionic fields $b(z)$ and $c(z)$ 
expressed in terms of the scalar
field $\ps$. Furthermore we 
can construct two particular fields as some
normally ordered polynomials of fields $b(z)$ and $c(z)$
and their derivative fields:
a Virasoro field $T(z)$ of conformal weight $2$
and a field $W(z)$ of conformal weight $3$.
These two fields $T(z)$ and $W(z)$ satisfy the operator
product expansion of the
$\Wth$ algebra with central charge $-2$.
The three fields $J^0 (z)$, $T(z)$ and $W(z)$
may be regarded as generating fields of
a $\WGth$ algebra. We will show that all the 
$J^i (z) =  \bIg, \; i = 0, 1, \ldots, $
can be expressed (See Lemmas \ref{lem_2} and \ref{lem_3})
as some normally ordered 
polynomials in terms of $T(z)$, $W(z)$ and
$J^0 (z)$ and their derivative fields. Since the space
$\Fbar \bigotimes \Hil$, being isomorphic to
$\Ml$, is an irreducible module
over the vertex algebra $\Winfone$, it is
also irreducible as a module over the $\WGth$ algebra. 

One can show 
that $J^0 (z) = i \partial \ph,$
by using Friedan-Martinec-Shenker bosonization.
Note that when the $\WGth$ algebra acts on
$\Fbar \bigotimes \Hil$, the Fourier components of
fields $T(z)$ and $W(z)$ act only on the first factor
$\Fbar $ while $J^0 (z)$ acts only on the second factor $\Hil$.
This implies that $\Fbar$ is irreducible as a module
over the $\Wthtwo$ algebra. We obtain full
character formulas of these irreducible modules $\Fbar$ of the
$\Wthtwo$ algebra
as a consequence of our explicit free field realization.
As a by-product of our free field
realization of $\Fbar$, we find that there
exists non-split short exact sequences
of modules over the $\Wth$ (resp. $\Winfone$) algebra with
central charge $-2$ (resp. $-1$).

The plan of this paper goes as follows. 
In Section \ref{sect_Woneinfty}, we review
the definition of $\hD$ and the construction
of vertex algebra $\Wc$. We present the free field
realization of $\Winfone$ in terms of $\bg$
fields. In Section \ref{sect_boson}, we 
recall the bosonization of $\bg$ fields in detail.
In Section \ref{sect_W3} we review the $\Wth$ algebra
in the framework of vertex algebras. 
In Section \ref{sect_relation} we
prove that vertex algebra $\Winfone$
is isomorphic to a tensor product
of the simple $\Wthtwo$ algebra and 
a Heisenberg vertex algebra 
generated by a free bosonic field.
We construct a number of irreducible modules of 
the $\Wthtwo$ algebra. In Section \ref{sect_char}, 
we calculate the full character
formula for representations of the $\Wthtwo$
algebra constructed in Section \ref{sect_relation}.

We take this opportunity to announce that we classify
the irreducible modules of the
$\Wthtwo$ algebra in our subsequent paper 
\cite{W2}. It turns out that these irreducible modules
are parametrized by points on a certain rational curve. We will
also classify all the irreducible modules of
$\Winfone$ algebra based on the relation between 
$\Wthtwo$ and $\Winfone$ algebras found in this paper.

\section{Vertex algebra $\Wc$ and 
  free fields realization of $\Winfone$}
\label{sect_Woneinfty}

Let $\D$ be the Lie algebra of regular differential operators on
the circle. The elements 
\begin{eqnarray*}
  J^l_k = - t^{l+k} ( \partial_t )^l, 
     \quad l \in \Z_{+}, k \in \Z, 
\end{eqnarray*}
form a basis of $\D$. $\D$ has also another basis
\begin{eqnarray*}
  L^l_k = - t^{k} D^l, 
     \quad l \in \Z_{+}, k \in \Z,  
\end{eqnarray*}
where $D = t \partial_t$. Denote by $\hD$ the central extension of  
$ \D $
by a one-dimensional center with a generator $C$, with
commutation relation (cf. \cite{KR1})
\begin{eqnarray}
  \left[
     t^r f(D), t^s g(D)
  \right]
    & = & t^{r+s} 
    \left(
      f(D + s) g(D) - f(D) g(D+r) 
    \right) \nonumber \\ 
   & + & \Psi 
      \left(
        t^r f(D), t^s g(D)
      \right)
      C,
  \label{eq_12}
\end{eqnarray}
where 
\begin{equation}
  \Psi 
      \left(
        t^r f(D), t^s g(D)
      \right)
   = 
   \left\{
      \everymath{\displaystyle}
      \begin{array}{ll}
        \sum_{-r \leq j \leq -1} f(j) g(j+r),& r= -s \geq 0  \\
        0, & r + s \neq 0. 
      \end{array}
    \right. \\
  \label{eq_13}
\end{equation}

Letting weight $J^l_k = k$ and weight $ C = 0$ defines a principal
gradation
\begin{equation}
  \hD = \bigoplus_{j \in \Z} \widehat{\cal D}_j.
  \label{eq_14}
\end{equation}
Then we have the triangular decomposition
\begin{equation}
  \hD = \hD_{+} \bigoplus \hD_{0}  \bigoplus \hD_{-},
  \label{eq_15}
\end{equation}
where 
\begin{eqnarray*}
  \hD_{\pm} = \bigoplus_{j \in \pm \Bbb N} \widehat{\cal D}_j,
  \quad
  \hD_{0} = {\D}_{0}  \bigoplus {\Bbb C} C.
\end{eqnarray*}

Let $\cal P$ be the distinguished parabolic subalgebra of $\D$,
consisting of the differential operators that extends
into the whole interior of the circle. 
${\cal P}$ has a basis $\{ J^l_k, l \geq 0, l + k \geq 0 \}.$
It is easy to check that the $2$-cocycle $\Psi$ defining
the central extension of $\hD$ vanishes when restricted
to the parabolic subalgebra $\cal P$. So $\cal P$ is also a
subalgebra of $\hD$. Denote
$\widehat{\cal P} = {\cal P} \oplus {\Bbb C} C$.

Fix $c \in \Bbb C$. Denote by ${\Bbb C}_c$ the $1$--dimensional
$\widehat{\cal P}$ module by letting $C$ acts as scalar $c$ and 
$\cal P$ acts trivially. Fix a non-zero vector
$v_0$ in ${\Bbb C}_c$. The induced $\hD$--module
\begin{eqnarray*}
  M_c \left( \hD
      \right)
  = {\cal U} \left( \hD
             \right)
    \bigotimes_{ {\cal U} \left( \cal P
                          \right)
               }
    {\Bbb C}_c
\end{eqnarray*}
is called the vacuum $\hD$--module with central charge $c$. Here we 
denote by ${\cal U} (\frak g)$ the universal enveloping
algebra of a Lie algebra $\frak g$. $M_c (\hD)$ admits a unique
irreducible quotient, denoted by $\Wc$. 
Denote the highest weight vector $1\otimes v_0$ in
$M_c (\hD)$ by $\vac$.

It is shown in 
\cite{FKRW} that $\Wc$ carries a canonical vertex
algebra structure, with vacuum vector $\vac$ and
generating fields 
\begin{eqnarray*}
  J^l (z) = \sum_{k \in \Bbb Z} J^l_k z^{-k-l-1},
\end{eqnarray*}
of conformal weight $l + 1, l = 0, 1, \cdots.$ The fields
$J^l (z)$ corresponds to the vector $J^l_{-l-1} \vac$
in $\Wc$. Below we will concentrate on the 
particular case $\Winfone$.
       
Recall that the bosonic $\bg$ fields are
\begin{equation}
 \beta (z) = \SUM \beta \N,
 \quad
 \gamma (z) = \SUM \gamma \NN 
 \quad (s \in \Bbb C)
  \label{bega}
\end{equation}
with the operator product expansions (OPEs)
\begin{equation}
  \beta (z) \gamma (w) \sim \frac{-1}{z-w}(\frac{z}{w} )^s, \;
   \beta (z) \beta (w) \sim 0, \;
    \gamma (z) \gamma (w) \sim 0.
  \label{bg_OPE}
\end{equation}
In other words, we have the following
commutation relations
$$ \left[
     \gamma (m), \beta (n)
   \right]
     = \delta_{m,-n},  \quad
   \left[
     \beta (m), \beta (n)
   \right]
     = 0,  \quad
   \left[
     \gamma (m), \gamma (n)
   \right]
     = 0.
$$

Let us denote by $\M$ the Fock space of the $\bg$ 
fields, with the vacuum vector $\vacs $, and 
\begin{equation}
  \beta (n+1) \vacs = 0, \quad
  \gamma (n) \vacs = 0, \quad n \geq 0.  
 \label{eq_vacuum}
\end{equation}
One can realize a representation of $\Winfone$
on $\M$ by letting (cf. \cite{KR2, M}, our convention here is
a little different):
\begin{equation}
  J^N (z) = \bNg + \frac{1}{N+1} s(s-1)\cdots (s-N)z^{-N-1}, 
  N \in \Z_{+}.
  \label{eq_field}
\end{equation}
The normal ordering $::$ is understood as moving
the operators annihilating $\vacs$ to the right.

Note that 
$J^0 (z) = \sum_{ k \in \Z} J^0_k z^{-k-1}$
is a free bosonic field of conformal weight $1$
with commutation relations
$$ [ J^0_m, J^0_n ] = -m \delta_{m -n}, \; m,n \in \Z.  $$
We also have the following commutation relations:
$$ [ J^0_m, \beta (n)] =  \beta (m+n ), \quad
   [ J^0_m, \gamma (n)] = - \gamma (m+n),   \quad m,n \in \Z.   $$
Then we have
the $\bg$-charge decomposition of $\M$ according to the
eigenvalues of the operator $- J^0_0$:
$\M = \bigoplus_{ l \in \Z } \Ml. $ It is known \cite{M, KR2} that 
${\cal M}_0^0$ is isomorphic to $\Winfone$ as vertex algebras.

\section{Bosonizations}
  \label{sect_boson}
In Section \ref{subsect_fermion} we recall 
the well-known boson-fermion correspondence  (cf. \cite{F1}). 
In Section \ref{subsect_boson} we review the 
Friedan-Martinec-Shenker bosonization of
the $\bg$ fields and some more detailed
structures \cite{FMS, FF}.

\subsection{Bosonization of fermions}
\label{subsect_fermion}

Let $j(z)$ be a free bosonic field of conformal weight $1$,
namely
$$ j(z) j(w) \sim \frac{1}{(z-w)^2}, $$
or equivalently, by introducing $j(z) = \SUM j(n) z^{-n-1}$, we have
$$[j(m), j(n)] = m \delta_{m, -n}.$$

Let us also introduce the free scalar field
$$\ph = q + j(0) \ln z - \sum_{n \neq 0} \frac{j(n)}{n} z^{-n},$$ 
where the operator $q$ satisfies
$[q, j(n)] = \delta_{n,0}.$ Clearly $j (z) = \partial \ph.$

Given $\alpha \in \Bbb C$, we denote by ${\cal H}_{\alpha}$
the Fock space of the free field $j(z)$ generated by the
vacuum vector $\vacuum $ satisfying 
$$j(n)  \vacuum
  = \alpha \delta_{n,0} \vacuum, 
   \quad n \geq 0.$$
It is well known that $\Ho$ is a vertex algebra, which
we refer to as {\em a Heisenberg vertex algebra}.
Easy to see that 
$\exp (\eta q ) \vacuum  = \mid \alpha + \eta \rangle$. 
Introduce the vertex operator
$$ X_{\eta} (z) 
    = \SUM exp(\eta q) z^{\eta \alpha} X_{\eta} (n) z^{-n}   $$ 
as follows. Let 
\begin{eqnarray}
X_{\eta} (z) & = & : \exp \left(
                        \eta \ph
                      \right) :      \\
             & = & \exp (\eta q) z^{\eta \alpha}
                 \exp \left(
                        \eta \sum_{n > 0} j (-n) z^n/n 
                      \right)
                 \exp \left(
                        \eta \sum_{n < 0} j (-n) z^n/n 
                      \right). \nonumber
\label{def_vertex}
\end{eqnarray}
The Fourier components of
$X_{\eta} (z)$ act from ${\cal H}_{\alpha}$ to
${\cal H}_{\alpha + \eta}$. Furthermore we have
the following OPE
\begin{equation}
  j(z) X_{\eta} (w) \sim 
    \frac {\eta X_{\eta} (w)}{z-w} 
      + \frac{1}{\eta} \partial X_{\eta} (w),
\label{def_ope}
\end{equation}
or equivalently we have
\begin{eqnarray}
     \left[
       j(m), X_{\eta} (n) 
     \right] & = & \eta X_{\eta} (m + n), \nonumber \\
     : j(z) X_{\eta} (z) : & = & 
        \frac{1}{\eta} \partial X_{\eta} (z).  \nonumber     
\end{eqnarray}
Also we have 
\begin{equation}
  X_{\xi} (z) X_{\eta} (w) 
 \sim (z-w)^{\xi \eta} : X_{\xi} (z) X_{\eta} (w) :
   \label{ope_vertex}
\end{equation}
In particular we have a pair of fermionic fields $X_{\pm} (z)$
with OPEs:
\begin{eqnarray*}
     X_1 (z) X_{-1} (w) \sim \frac{1}{z-w}, \quad
     X_{\pm 1} (z) X_{\pm 1} (w) \sim 0. 
\end{eqnarray*}
This is the well-known boson-fermion correspondence.

\subsection{Bosonization of bosons}
\label{subsect_boson}

First let us introduce the $bc$ fermionic fields 
$$ b(z) = \SUM b (n) z^{-n}, \quad  c(z) = \SUM c (n) z^{-n-1} $$
with OPEs
\begin{equation}
b(z) c(w) \sim \frac{1}{z-w},
\quad b(z) b(w) \sim 0,
\quad c(z) c(w) \sim 0.
   \label{ope_bc}
\end{equation}
In other words, we have 
$$ \left[ b(m), c(n) 
   \right]_{+} 
     = \delta_{m, -n}, \quad
   \left[ b(m), b(n)
   \right]_{+} = 0, \quad
   \left[ c(m), c(n)
   \right]_{+} =0.
$$
We denote by $\F$ the Fock space of the $bc$ fields, generated by
$\vacbc$, satisfying

$$ b( n+1 ) \vacbc = 0, 
   \quad  c(n) \vacbc = 0, \quad n \geq 0. $$

Then 
$$ j^{bc} (z) = \/ : c(z) b(z) :\/ = \SUM j^{bc}_n z^{-n-1} $$ 
is a free boson of conformal weight $1$ with commutation relations
$$[j^{bc}_m, j^{bc}_n] = m \delta_{m, -n}, \quad m,n \in \Z. $$ 
We further have the following commutation relations:
$$ [j^{bc}_m, b(n)] = -b(m+n), \quad
   [j^{bc}_m, c(n)] =  c(m+n), \quad m,n \in \Z.
$$
Then we have the $bc$--charge decomposition of $\F$ 
according to the eigenvalues of $j^{bc}_0$:
$$ \F = \bigoplus_{l \in \Z} \Fl .$$
Following \cite{FF}, we consider the vector space 
\begin{eqnarray*}
 N(s) = \sum_{ l \in \Z} \Fl \bigotimes 
            {\cal H}_{i(s+l)},
\end{eqnarray*}
and we define the actions of $\beta (n), \gamma (n), n \in \Z$ 
on $N(s)$ by letting \cite{FMS}
\begin{eqnarray}
    \beta (z) & = & \SUM \beta \N = {\partial} b(z) X_{-i} (z), 
       \label{eq_beta}  \\
    \gamma (z) & = & \SUM \gamma \NN = c(z) X_{i} (z).
       \label{eq_gamma}  
\end{eqnarray}
It can be easily shown that the bosonic fields
$\beta (z), \gamma (z)$ defined above indeed satisfy the 
OPEs (\ref{bg_OPE}). The vector 
$\vacbc \bigotimes | is \rangle$
satisfies the vacuum condition (\ref{eq_vacuum}) by means of
(\ref{eq_beta}) and (\ref{eq_gamma}).
Then we have a homomorphism 
${\epsilon}: \M \longrightarrow N(s) $
as modules of the Heisenberg algebra spanned by
$\beta (n), \gamma(n), n \in \Z, $
by letting 
$$ {\epsilon} \left(
            \vacs
          \right)
  = \vacbc \bigotimes | is \rangle.  $$

This homomorphism is obviously an embedding since 
$\M$ is an irreducible module of the
above Heisenberg algebra.
The following proposition (cf. \cite{FF}) tells us the precise
image of this embedding. we reproduce the proof here since 
some crucial misprints in their proof in \cite{FF} 
need to be corrected.

\begin{proposition}
    The image of the homomorphism $\epsilon$ coincides with the
    kernel of $c(0)$, acting from $N(s)$ to $N(s-1)$.
  \label{prop_sub}
\end{proposition}

\begin{demo}{Proof}
  The operators $\beta (n), \gamma(n), \/ n \in \Z, $ 
given by (\ref{eq_beta}) and (\ref{eq_gamma}) do not depend on 
$b(0)$ and therefore commute with $c(0)$. 
So we have $ \mbox{Im} \/ \epsilon \subset \ker c(0)$
since the operator $c(0)$ kills the vacuum vactor
$ \vacbc$. It is easy to
see that the kernel of $c(0)$ is obtained by applying to
$ \vacbc \bigotimes \mid is \rangle $
the operators $j(n), c(n), n \in \Z,$ and 
$ b(m), m \in \Z - \{ 0 \}.$
So it remains to show that fields $j(z), c(z)$ and $\partial b(z)$ 
can be expressed in terms of fields $\beta (z), \gamma (z)$,
$X_{\pm i} (z)$ and their derivative fields. 
Indeed it is easy to show that 
\begin{eqnarray*}
    \partial b(z) & = & \partial \beta (z) X_i (z), \\
    c(z) & = & \partial \gamma (z) X_{-i} (z).
\end{eqnarray*}
Recall that $J^0 (z) = :\gamma (z) \beta (z):$. Easy to
check by (\ref{eq_beta}) and (\ref{eq_gamma}) that
$j (z) \equiv \partial_z \ph = - i J^0 (z)$.
\end{demo}

Denote by $ \Fbar $ the kernel of the operator $c(0)$
acting from $ \Fl $ to ${\cal F}^{l+1} $.
We now have a natural isomorphism:
\begin{equation}
  \Ml \cong \overline{\cal F}^l \bigotimes \Hil.
 \label{eq_isom}
\end{equation}

\begin{remark}
  $\overline{\cal F}^0$ is a vertex subalgebra of $ {\cal F}^0$.
 This is an example of the following well-known fact in
 the theory of vertex algebras: Given a vertex algebra
 $V$ and let $Y(a, z) = \SUM a(n) z^{-n-1}$ be 
 the field corresponding  to some vector $a \in V$, 
 then the kernel of the operator $a(0)$ acting on $V$
 is always a vertex subalgebra of $V$.
\end{remark}

\section{$\Wth\/$ algebra}
  \label{sect_W3}

Denote by $\/\UW$ ($ c \in \Bbb C$ is the central charge)
the quotient of the
free associative algebra 
generated by $L_m, W_m,\; m \in \Z$ by the two-sided ideal 
generated by the following commutation relations 
(cf. e.g. \cite{BMP}):
\begin{eqnarray}
{[ L_m, L_n ]} & = &  
    (m-n) L_{m+n} + \frac{c}{12}(m^3 - m) \delta_{m, -n}, \nonumber  \\
{[ L_m, W_n ]} & = & 
     (2m - n) W_{m + n},                \nonumber      \\
{[ W_m, W_n ]} & = &    
    (m-n)\Bigl(
           \frac{1}{15}(m+n+3)(m+n+2)     \label{eq_alg}      \\
     && \quad\quad\quad\quad\quad\quad\quad\quad\quad - \frac{1}{6} (m+2)(n+2)
         \Bigl) L_{m+n}                   \nonumber         \\
     && + \beta (m-n) \Lambda_{m+n} 
        + \frac{c}{360}m(m^2 -1)(m^2 -4) \delta_{m,-n}, \nonumber    
\end{eqnarray}
with $\beta = 16/(22+5c) $ and
\begin{eqnarray*}
   \Lambda_m = \sum_{k \leq -2} L_k L_{m-k} 
               + \sum_{k > -2} L_{m-k} L_k 
               - \frac{3}{10} (m+2)(m+3) L_m.
\end{eqnarray*}

Denote
\begin{eqnarray*}
   \Wthpm = \{ L_n, W_n, \; \pm n \geq 0 \}, \quad
   \Wtho = \{ L_0, W_0 \}.
\end{eqnarray*}

A Verma module $\/ \Mtw$ of $\/\UW\/$ is the
induced module
\begin{eqnarray*}
  \Mtw = 
   \UW \bigotimes_{ {\cal U}( \Wthp \oplus \Wtho) }
    {\Bbb C}_{t, w}
\end{eqnarray*}
where $\/{\Bbb C}_{t, w}$ is the 1-dimensional module of
$\/{\cal U}( \Wthp \oplus \Wtho)$ such that 
\begin{equation}
  \Wthp \vactw = 0, \,L_0 \vactw = t \vactw, \,
   W_0 \vactw = w \vactw.
\end{equation}
$\Mtw $ has a unique irreducible quotient which is denoted by
$\Ltw$. A singular vector in a $\/ \UW$-module
means a vector killed by $\Wthp$. It is easy to see that
$L_{-1} \vacoo, W_{-1} \vacoo$, and $W_{-2} \vacoo$
are singular vectors in $\Moo$. We denote
by $\VM$ the {\em vacuum module} which is by definition the quotient
of the Verma module $\Moo\/$ by the $\/ \UW$-submodule 
generated by the singular vectors 
$L_{-1} \vacoo, W_{-1} \vacoo$, and $W_{-2} \vacoo$.
We call ${\cal L}_c (0, 0)$ the {\em irreducible vacuum module}.
Let $I$ be the maximal proper submodule of 
the vacuum module $\VM$.
Clearly ${\cal L}_c(0, 0)$ is the irreducible quotient
$\VM / I$ of $\VM$. It is easy to see that $\VM$ has a linear basis 
\begin{eqnarray}
  L_{-i_1 -2} \cdots L_{-i_m -2} W_{-j_1 -3} \cdots W_{-j_n -3} \vacoo, 
                    \nonumber   \\
    0 \leq i_1 \leq \cdots \leq i_m, 
   \quad 0 \leq j_1 \leq \cdots \leq j_n, \quad m,n \geq 0. 
\end{eqnarray}

Introduce the following fields
\begin{equation}
  T(z) = \sum_{n \in \Z} L_n z^{-n-2},  \quad
  W(z) = \sum_{n \in \Z} W_n z^{-n-3}.
\end{equation}
It is well known that the vacuum module
$\VM$ (resp. irreducible vacuum module ${\cal L}_c(0, 0)$)
carries a vertex algebra structure 
with generating fields $T(z)$ and $W(z)$.
The $\Wth$ algebra with central charge $-2$ we have
been referring to
is the vertex algebra ${\cal L}_{-2} (0, 0)$, which we denote
by $\Wthtwo$ throughout our paper. 
Fields $T(z)$ and $W(z)$ correspond to the 
vectors $L_{-2} \vacoo$ and $W_{-3} \vacoo$ 
respectively. The field corresponding to the vector
$L_{-i_1 -2} \cdots L_{-i_m -2} W_{-j_1 -3} \cdots W_{-j_n -3} \vacoo$
is 
$$ \/\partial^{(i_1)} T(z) \cdots \partial^{(i_m)} T(z) 
  \/\partial^{(j_1)} W(z) \cdots \partial^{(j_n)} W(z), $$
where $\/\partial^{(i) }\/$ denotes 
$ \/\frac{1}{i!}\partial^i_z $. We can rewrite (\ref{eq_alg})
in terms of the following OPEs in our central charge $ -2$ case:
  \begin{eqnarray}
    T(z) T(w) & \sim & \frac{-1}{(z-w)^4} 
                 + \frac{2 T(w)}{ (z-w)^2 } 
                 + \frac{\partial T(w)}{z-w}    \nonumber \\
    T(z) W(w) & \sim & \frac{3W(w)}{(z-w)^2}
                     + \frac{\partial W(w)}{z-w} \label{ope_three}\\
    W(z) W(w) & \sim & \frac{-2/3}{(z-w)^6}
                     + \frac{2T(w)}{(z-w)^4}
                     + \frac{\partial T(w)}{(z-w)^3} \nonumber \\
                 &&  + \frac{1}{(z-w)^2}
                        \left(
                         \frac{8}{3}:T(w)T(w): - \hf {\partial}^2 T (w)
                        \right) \nonumber \\
                 &&  + \frac{1}{z-w} 
                        \left(
                           \frac{4}{3} \partial \left(
                                               :T(w) T(w):
                                               \right)
                           - \frac{1}{3} {\partial}^3 T (w)
                        \right).                    \nonumber 
  \end{eqnarray}

Representation theory of the vertex algebra $\VM$ is
just the same as that of $\/ {\cal U}(\Wth)$.
However note that $c = -2$ is not a generic central
charge \cite{W2}, namely the vacuum module $\VM$ with $c = -2$
is reducible, or in other word, the maximal proper submodule
$I$ of $\VM$ is not zero. 
Thus representation theory of $\Wthtwo$ becomes
highly non-trivial due to the following constraints:
a module $M\/$ of the vertex algebra
$\VM$ can be a module of the vertex algebra
$\Wthtwo$ if and only if $M\/$ is annihilated
by all the Fourier 
components of all fields corresponding to vectors
in the maximal proper submodule $I \subset \VM$.

\section{Relations between $\Wth$ algebra
and vertex algebra $\Winfone$}
   \label{sect_relation}

Define 
\begin{equation}
  T (z) \equiv \SUM L_n z^{-n-2}
    = : \partial b (z) c(z) :. 
  \label{eq_vir}
\end{equation}
Easy to check that 
$T(z)$ is a Virasoro field with central charge $-2$.
We also define another field of conformal weight $3$:
\begin{equation}
  W(z) \equiv \SUM W_n z^{-n-3}
         = \frac{1}{\sqrt{6}}
              \left(
                : {\partial}^2 b(z) c(z) : 
                  - :\partial b(z) \partial c(z) :
              \right). 
   \label{eq_w3}
\end{equation}

We have the following proposition whose proof is straightforward
however tedious by using Wick's theorem.

\begin{proposition}
   Fields $ T(z)\/$ and $W(z)\/$ satisfy the OPEs (\ref{ope_three})
  of $\Wth\/$ algebra with central charge $-2$.
 \label{prop_wthr}
\end{proposition}

We note that this $\Wth$ algebra structure in $bc$ fields
was also observed in \cite{BCMN}. 

We rescale $W(z)$ to be 
$\widetilde{W} (z) = \hf \sqrt{6} \/ W(z)$, namely
\begin{equation}
 \widetilde{W} (z)  \equiv \SUM 
  \widetilde{W}_n z^{-n-3}
   = \hf \left(
                : {\partial}^2 b(z) c(z) : 
                  - :\partial b(z) \partial c(z) :
         \right).
  \label{eq_wtilde}
\end{equation}
We will see later that it is more convenient to work with the
rescaled field $\widetilde{W} (z)$.
Now we can state our first main results of this paper.

\begin{theorem}
   \begin{enumerate}
 \item[1)] The vertex algebra $\overline{\cal F}^0$ is isomorphic
to the simple vertex algebra $\Wthtwo$, with generating
fields $T(z)$ and $W(z)$. $\Fbar$ $(l \in \Z)$ are irreducible
modules of the $\Wthtwo$ algebra.
 
  \item[2)] The vertex algebra $\Winfone$ is isomorphic to a tensor
product of the vertex algebra $\Wthtwo$, and the Heisenberg 
vertex algebra ${\cal H}_0$ with $ J^0 (z)$ as 
a generating field. 
\end{enumerate}
  \label{th_theorem}
\end{theorem}

Proof of the above theorem reply on the following three
lemmas:

\begin{lemma}

   The vector space $\Ml$
is irreducible regarded as a module of the vertex algebra
$\Winfone$ via the free field realization (\ref{eq_field}).
  \label{lem_1}
\end{lemma}

Lemma \ref{lem_1} was proved in \cite{KR2, M}.

\begin{lemma}
  The fields $J^n (w) = \bngw, n \geq 0$ acting on the Fock space
${\cal M}_0$ can be expressed as a normally
ordered polynomial in terms of fields
$$ :{\partial}^i b(w) {\partial}^j c(w)~:~, i + j \leq n, i > 0,
   \mbox{ and } 
   {\partial}^k j(w), k = 0, 1, \cdots, n.  $$
More precisely, we have
  \begin{equation}
    \begin{array}{rcl}
      \lefteqn{
        \bngw  }       \\
      & = & \sum_{ 1 \leq k \leq n}
              \left[
                k \bi :\partial^{n-k+1} b(w) c(w):P_{k-1}(j)
              \right]
          + C_n P_{n+1} (j),
    \label{eq_key}
    \end{array}
  \end{equation}
where $C_n = (n+2) \sum_{m=0}^n (-1)^{m+1}\frac{1}{m+2}$ is 
some constant depending on $n$, and
the normally ordered polynomial $P_m (j)$ (or denoted by
$P_m ( j(w) )$ when it is necessary to specify
the variable in $j(w)$) in terms of
the field $j(w)$ and its derivative fields is defined as (recall 
that $j(w) = \partial \phi (w)$)
  \begin{equation}
    P_m (j) = \frac { \partial_w^m :e^{-i \phi (w) }: }{:e^ {- i \phi (w)} :},
    \quad m \geq 0.
   \label{eq_der}
  \end{equation}
 \label{lem_2}
\end{lemma}

\begin{demo}{Proof of Lemma \ref{lem_2}}

  We will calculate the normally ordered product
$ : {\partial}^n \beta(w) \gamma (w) : $ instead of $\bngw$.
These two normally ordered products coincide since both
$\beta (w)$ and $ \gamma (w)$ are free fields.

By formulas (\ref{eq_beta}) and (\ref{eq_gamma}), we have
  \begin{eqnarray}
    : {\partial}^n \beta(w) \gamma (w) :
      & = & : {\partial}^n \left(
                         \partial b(w) X_{-i} (w) 
                           \right)
           (c(w) X_i (w)) :              \label{eq_fms}        \\
      & = & \sum_{ 0 \leq k \leq n} \bi 
              \left[
               : {\partial}^{n-k+1} b(w) {\partial}^k X_{-i} (w)
               c(w) X_i (w) :
              \right].\nonumber 
   \label{eq_36}
  \end{eqnarray}

It follows from the OPEs (\ref{ope_bc}) that
  \begin{eqnarray}
    {\partial}_z^{n-k+1} b(z)  c(w) & = &
      \frac{(-1)^{n-k+1} (n-k+1)!} {(z-w)^{n-k+2}}
       + :{\partial}^{n-k+1} b(w) c(w): \nonumber  \\
       & & +\mbox{\,\,higher terms},   
   \label{eq_37} 
  \end{eqnarray}

It follows from the OPEs (\ref{ope_vertex}) that
  \begin{eqnarray}
    {\partial}_z^k X_{-i} (z) X_{i} (w)
     & = &
      {\partial}^k_z \left(
                       \sum_{m \geq 0} \frac{(z-w)^{m+1}}{m!} P_m (j(w))
                     \right)               \label{eq_38}    \\
     & = & \sum_{m \geq k-1}\frac{m+1}{(m-k+1)!} 
             (z-w)^{m-k+1} P_m (j(w) ).\nonumber    
  \end{eqnarray}

Since $ : {\partial}^n \beta(w) \gamma (w) : $ 
is the constant term in the expansion of 
power series of $z-w$ in the
operator product expansion of 
${\partial}^n \beta(z) \gamma (w)$, we see
from equations (\ref{eq_fms}), (\ref{eq_37}) and 
(\ref{eq_38}) that the only terms in
equation~(\ref{eq_38}) which will contribute to 
$ : {\partial}^n \beta(w) \gamma (w) : $ non-tivially
is the two terms $m = k-1$ and $m = n+1$.
Namely we have
  \begin{equation}
    \begin{array}{rcl}
      \lefteqn{
        : {\partial}^n \beta(w) \gamma (w) : 
              }           \\
       &= &\sum_{ 0 \leq k \leq n} \bi 
              \Bigl[
                k : {\partial}^{n-k+1} b(w) c(w) : P_{k-1} (j) \\
       & &    + \frac{n+2}{(n-k+2)!}
                 (-1)^{n-k+1} (n-k+1)! P_{n+1} (j)
              \Bigl]                                           \\
       &= & \sum_{ 1 \leq k \leq n} 
              \left[
                k \bi : {\partial}^{n-k+1} b(w) c(w) : P_{k-1} (j)
              \right]
             + C_n P_{n+1} (j), \nonumber
    \end{array}
  \end{equation}
where 
$$ C_n = (n+2) \sum_{m=0}^n (-1)^{m+1}\frac{1}{m+2}.  $$
\end{demo}

\begin{remark}
 \begin{enumerate}
  \item[1)] $P_m (j)$ defined in equation (\ref{eq_der})
   reads as follows for small $m$:
   \begin{eqnarray*}
    P_1 (j) & = & -i j(w), \quad P_0 (j) = 1, \\
    P_2 (j) & = & -i \partial j(w) - :j(w)^2:,\\
    P_3 (j) & = & -i \partial^2 j(w) -3 :j(w) \partial j(w): + i :j(w)^3:.
   \end{eqnarray*} 

  \item[2)] The formula (\ref{eq_key}) reads as follows for small $n$:
  \begin{eqnarray*}
  : \gamma (w) \beta (w) : & \equiv & J^0 (w) = i j(w),  \nonumber      \\
  : \gamma (w) \partial \beta (w) : & = & 
     : \partial b (w) c(w) : - \hf :J^0 (w)^2 : 
       + \hf \partial J^0 (w),       \\
  : \gamma (w) \partial^2 \beta (w) : 
     & = & 2 : \partial^2 b(w) c(w): - 2 : \partial b(w) c(w) : J^0 (w) \\
       & & + \frac{5}{3} \partial^2 J^0 (w) 
        + \frac{5}{3} : J^0 (w)^3 :        \\ 
       & & - 5 :J^0 (w) \partial J^0 (w):.
  \end{eqnarray*}
 \end{enumerate}
\end{remark}
\begin{lemma}
   Each field $:{\partial}^i b(z) {\partial}^j c(z) :, i > 0, j \geq 0$
can be expressed as a normally ordered polynomial 
in terms of $\,T(z)$ and $\/W(z)$ defined in (\ref{eq_vir})
and (\ref{eq_w3}) and their derivative fields.
  \label{lem_3}
\end{lemma}

\begin{demo}{Proof of Lemma \ref{lem_3}}
  We first prove the following statement:

Claim $A_n$: Any field 
$:{\partial}^i b(z) {\partial}^{n-i+1} c(z):, 1 \leq i \leq n+1$
can be written as a linear combination of the following
$n+1$ fields
$$ \partial \left(
            :{\partial}^k b(z) {\partial}^{n-k} c(z):
          \right),
  1 \leq k \leq n \mbox{  and }
 :T (z) {\partial}^{n-1} b(z) c(z) :. $$

Indeed one can calculate directly by using (\ref{eq_vir})
and Wick's Theorem that
\begin{equation}
:T (z) \left( 
          {\partial}^{n-1} b(z) c(z)
       \right) :
  = \hf :  {\partial}^{n-1} b(z) {\partial}^2 c(z) :
    + \frac{1}{n} :{\partial}^{n+1} b(z) c(z) :.
\label{eq_31}
\end{equation}
Also since the derivation of a normally ordered product
satisfies the Leibniz rule we have
\begin{equation}
 \partial \left(
             : {\partial}^{k} b(z) {\partial}^{n-k} c(z) :
           \right)
  = : {\partial}^{k+1} b(z) {\partial}^{n-k} c(z) :
    + : {\partial}^{k} b(z) {\partial}^{n-k+1} c(z) :.
\label{eq_32}
\end{equation}
The $n+1$ fields 
$$  \partial \left(
             : {\partial}^{k} b(z) {\partial}^{n-k} c(z) :
          \right),
    1 \leq k \leq n
  \mbox{ and } 
:T (z) {\partial}^{n-1} b(z) c(z) :  $$
can be obtained from the $n+1$ fields
$:{\partial}^i b(z) {\partial}^{n-i+1} c(z):, 1 \leq i \leq n+1$
through a linear transformation given by the following 
$(n+1) \times (n+1) $ matrix
\begin{eqnarray*}
\left[ \begin{array}{cccccc}
1      & 1      &        &             &        &       \\
0      & 1      & \ddots &             &        &       \\
       & \ddots & \ddots & \ddots      &        &       \\
       &        & \ddots & 1           & 1      & 0     \\
       &        &        & 0           & 1      & 1      \\
       &        &        & \frac{1}{2} & 0      & \frac{1}{n+1} 
\end{array} \right].
\end{eqnarray*}
It is easy to see this matrix has determinant $\frac{n+3}{2(n+1)}$
so it is invertible. By inverting
the matrix we prove the Claim $A_n$. 

Now we are ready to prove the following claim by induction
on $n$ which is a reformulation of Lemma \ref{lem_3}.
Claim $B_n$: Any field 
$:{\partial}^i b(z) {\partial}^{n-i} c(z), 1 \leq i \leq n$
can be written as a normally ordered polynomial
in terms of $T (z), W(z)$ and their derivative fields.

When $n = 1$, $: \partial b(z) c(z) : $ 
is just $ T(z) $ itself.

When $n = 2$, $: {\partial}^2 b(z) c(z) : $ and
$ : \partial b(z) \partial c(z) : $ are clearly
linear combinations of the fields 
$$ \partial T(z) = \partial \left(
            : \partial b(z) c(z) :
          \right)
  = : {\partial}^2 b(z) c(z) : 
    + : \partial b(z) \partial c(z) : $$
and 
$$ W(z) = \frac{1}{\sqrt{6} } \left(
                    : {\partial}^2 b(z) c(z) : 
                    - : \partial b(z) \partial c(z) :
                 \right). $$
So Claim $B_2$ is true.

Assume that the statement $B_n$ is true. 
Then particularly the field $: {\partial}^{n-1} b(z) c(z) :$
can be written as a normally ordered polynomial
of $ T (z)\/$ and $W(z)\/$. And so is 
$ :T(z) {\partial}^{n-1} b(z) c(z) : $.
Then the Claim $B_{n+1}$ follows from 
Claim $A_n$ (cf. equation (\ref{eq_31})).
\end{demo}

\begin{demo}{Proof of Theorem \ref{th_theorem}}
  Lemmas \ref{lem_1} and \ref{lem_2} imply immediately
  that $\Fbar$ is irreducible under the actions of 
  the Fourier components of fields
  $:{\partial}^i b(z) {\partial}^j c(z) :, i > 0, j \geq 0$.
  Together with Lemma \ref{lem_3}, this implies 
  that $\Fbar$ is irreducible under the actions of
  $L_n, W_n, \/ n \in \Z.$ So the vertex algebra
  $\overline{\cal F}^0$ is isomorphic to $\Wthtwo$ by 
  Proposition \ref{prop_wthr}. The free field construction of $\Fbar$ 
  guarantees that $\Fbar$
  is a module of the vertex algebra $\Wthtwo$.
  The second statement of Theorem \ref{th_theorem}
  now follows from the isomorphism of vertex algebras
  $M^0_0 \cong {\overline{\cal F}}^0 
  \bigotimes {\cal H}_0$ given by (\ref{eq_isom}).
\end{demo}

We have the following proposition from the explicit 
free field realization of modules $\Fbar$
of the $\Wthtwo$ algebra. Also see Remark~{4.3} in
\cite{W1} for some further implication.
\begin{proposition}
 \begin{enumerate}
 \item[1)] There exists a non-split short exact sequence 
  of modules over the vertex algebra $\Wthtwo$:
  \begin{equation}
    0 \longrightarrow \Fbar \longrightarrow \Fl 
       \longrightarrow {\overline{\cal F}}^{l-1} \longrightarrow 0. 
    \label{eq_split1}
  \end{equation}
 \item[2)] There exists a non-split short exact sequence 
  of modules over the vertex algebra $\Winfone$:
  \begin{equation}
    0 \longrightarrow {\cal M}^l_{s-l} 
       \longrightarrow {\cal M}
       \longrightarrow {\cal M}^{l-1}_{s-l+1}
       \longrightarrow 0.
    \label{eq_split2}
  \end{equation}
  Here ${\cal M}$ is isomorphic to
    $\Fl \bigotimes {\cal H}_{is}$ as vector spaces.
 \end{enumerate}
\end{proposition}

\begin{demo}{Proof}
  As a vector space we have a direct sum 
  ${\cal F}^l = \Fbar \bigoplus b(0) {\overline{\cal F}}^{l-1}$.
  Then it is not hard to see that as a $\Wthtwo$-module, 
  $\Fl/\Fbar$ is
  isomorphic to the irreducible module 
  ${\overline{\cal F}}^{l-1}$. So the following non-split 
  short exact sequence of modules over the $\Wthtwo$ algebra
  $$ 0 \longrightarrow \Fbar \longrightarrow \Fl \longrightarrow 
       \Fl/\Fbar \longrightarrow 0  $$
  is isomorphic to the one in (\ref{eq_split1}).
  
  Note that ${\cal M}_s^j $ is isomorphic to 
  $\overline{\cal F}^j \bigotimes {\cal H}_{i(s+j)}$ as 
  modules over the vertex algebra $\Winfone$ by Theorem \ref{th_theorem}.
  Then the non-split short exact sequence (\ref{eq_split2})
  can be obtained by tensoring the one in (\ref{eq_split1}) with
  ${\cal H}_{is}$.
\end{demo}

\section{Character formulas of modules over $\Wth$ algebra
with central charge $-2$}
   \label{sect_char}

Denote by 
$$ \Psi (z, q, p) \equiv \sum_{ l \in \Z} z^l \psi_l (q,p)
= \Tr \mid_{\bigoplus_{l \in \Z} \Fbar} 
    z^{-j^{bc}_0} q^{L_0} p^{\widetilde{W}_0}$$
the full character of 
$\bigoplus_{l \in \Z} \Fbar$, a direct sum
of irreducible modules $\Fbar$ over the $\Wthtwo$ algebra. Here
$\psi_l(q,p)$ is the full character of $\Fbar, l \in \Z$.
  Then the full character formula $\psi_l (q,p)$ of the irreducible
  $\Wthtwo$-module $\Fbar$ can be recovered 
  from $ \Psi (z, q, p)$ by taking residue
$$\psi_l(q,p) = \mbox{Res}_{z = 0} z^{l+1}\Psi (z, q, p).$$

We will need the following lemma.
\begin{lemma}
  We have the following OPEs:
     \begin{eqnarray}
       T(z) b(w) & \sim & \frac{ \partial b(w) }{ z-w }, \nonumber \\
       T(z) c(w) & \sim & \frac{ c(w) }{ (z-w)^2 } 
                     + \frac{ \partial c(w) }{ z-w },  \nonumber  \\ 
       \widetilde{W}(z) b(w) 
                 & \sim & \frac{\hf \partial b(w) }{ (z-w)^2 } 
                     + \frac{ {\partial}^2 b(w) }{ z-w }, \label{ope_Wb} \\
       \widetilde{W}(z) c(w) 
                 & \sim & \frac{ - c(w) }{ (z-w)^3 }
                     + \frac{ - \frac{3}{2} \partial c(w) }{ (z-w)^2 } 
                     + \frac{ - {\partial}^2 c(w) }{ z-w }. \label{ope_Wc}
  \end{eqnarray}
    \label{lem_41}
\end{lemma}

\begin{demo}{Proof}
   We will prove the OPE (\ref{ope_Wc}) only and the other OPEs can
be proved similarly by using Wick's Theorem.

Since 
\begin{eqnarray*}
  b(z) c(w) \sim \frac{1}{z-w}, 
\end{eqnarray*}
we have 
\begin{eqnarray*}
  \left(
    \partial^2_z b(z) \right) c(w)
   \sim \frac{2}{(z-w)^3}. 
\end{eqnarray*}

Since $ c(z)c(w) \sim 0 $ and
$b(z), c(z)$ are fermionic fields, we have by Wick's Theorem
\begin{eqnarray}
\left(
  \partial^2_z b(z) c(z) \right) 
    c(w) 
  &\sim & - \frac{2 c(z)}{(z-w)^3} \nonumber  \\
  &\sim & - \frac{2 c(w)}{(z-w)^3} - \frac{2 \partial c(w)}{(z-w)^2}
        - \frac{ \partial^2 c(w) }{z-w}. \label{ope_Wc1}
\end{eqnarray}

We also have by Wick's Theorem
\begin{eqnarray}
\left(
  \partial_z b(z) \partial_z c(z) \right)   c(w)  
     \sim  \frac{ \partial_z c(z)}{(z-w)^2}  
     \sim  \frac{ \partial_w c(w)}{(z-w)^2}
          + \frac{ \partial_w^2 c(w)}{z-w}. \label{ope_Wc2}
\end{eqnarray}

Now the OPE (\ref{ope_Wc}) follows
from (\ref{ope_Wc1}), (\ref{ope_Wc2})
and the definition of $\widetilde{W}(z)$ in (\ref{eq_wtilde}).
\end{demo}

In particular Lemma \ref{lem_41} implies

\begin{corollary}
  We have the following commutation relations ($n \in \Z $):
  \begin{eqnarray*}
    {[ L_0, b(n) ]} = -n b(n), \quad
    {[ L_0, c(n) }] = -n c(n),  \\    
    \quad \quad   {[ W_0, b(n) ]} = n^2 b(n), \quad \; 
    {[ W_0, c(n) ]} = -n^2 c(n).
  \end{eqnarray*}
    \label{cor_41}
\end{corollary}    

\begin{demo}{Proof}
  By comparing the coefficients of the $z^{-3}$ terms
  in both sides of the OPE (\ref{ope_Wb}), we get
  \begin{equation}
   [ W_0, b(w) ] = w \partial b(w) + w^2 \partial^2 b(w).
    \label{eq_comm}
  \end{equation}
  Comparing the coefficients of the $w^n$ terms in
  both sides of (\ref{eq_comm}), we get 
  $       [ W_0, b(n) ] =  n^2 b(n). $ 
  Proofs of the other commutation relations in Corollary 
  \ref{cor_41} are similar.
\end{demo}

The following full character formula follows now
from Corollary \ref{cor_41} and the characterization of 
$\Fbar$ as the subspace of $\Fl$ consisting of
vectors which do not involve $b(0)$, 
the zero-th Fourier component of the field $b(z)$.
\begin{theorem}
  The full character formula $\Psi (z, q, p)$ is given by
    \begin{eqnarray*}
      \Psi (z, q, p) =
            \prod_{ n \geq 1} 
                  \left( 1+ z      q^n p^{n^2} 
                  \right)
                  \left( 1+ z^{-1} q^n p^{- n^2} 
                  \right).
     \end{eqnarray*}
   \label{th_char}
\end{theorem}
\begin{remark}
 \begin{enumerate}
  \item[1)] By using the Jacobi triple identity, one can easily show 
  that
  \begin{eqnarray*}
     \psi_l (q, 0) = \frac{1}{\Pi_{n \geq 1} (1-q^n)}
     \sum_{k \geq |l|} (-1)^{k+l} 
       q^{ \frac{k(k+1)}{2} }.
   \end{eqnarray*}
  This is consistent with the explicit decomposition of $\Fbar$ with
  respect to the Virasoro algebra generated by the
  Fourier components of the field $T(z)$ \cite{FF}.

  \item[2)] If we consider instead $$ \widetilde{\Psi} (z, q, p) 
   \equiv \Tr \mid_{\bigoplus_{l \in \Z} \Fl} 
    z^{-j^{bc}_0} q^{L_0} p^{\widetilde{W}_0},  $$
 then we can show similarly that 
    \begin{equation}
      \widetilde{\Psi} (z, q, p) =
            (1+z) \prod_{ n \geq 1} 
                  \left( 1+ z      q^n p^{n^2} 
                  \right)
                  \left( 1+ z^{-1} q^n p^{- n^2} 
                  \right).
        \label{eq_modular}
     \end{equation}
  Essentially the same formula as in (\ref{eq_modular}) 
  up to some simple changes of variables
  appears in \cite{Di} as some
  generating function of counting 
  covers of an elliptic curve. 
  Modular invariance and some other interesting properties 
  of the function $\widetilde{\Psi} (z, q, p)$ were discussed
  in detail in \cite{KZ}. It is suggested that
  $\widetilde{\Psi} (z, q, p)$ may be an indication of the
  existence of generalized Jacobi forms involving several 
  (possibly infinitely many) variables.
  We hope that full character formulas
  of representations of $\cal W$-algebras in general
  may provide further natural examples of generalized Jacobi forms.
 \end{enumerate}
\end{remark}
{\bf Acknowledgement} The results of
this paper were presented in the Seminar of
Geometry, Symmetry and Physics at Yale University
and in the 1997 AMS Meeting at Detroit.
I thank the organizers of the meeting, Chongying Dong and Bob Griess
for invitation. I thank Edward Frenkel,
Igor Frenkel, Victor Kac, and
Gregg Zuckerman for their interests and comments, and especially
Edward Frenkel for stimulating discussions. 
I also thank Gerd Mersmann for pointing out to me
the references \cite{Di, KZ}.

\end{document}